\documentclass[a4paper, reqno]{amsart}

\usepackage{amsmath,amssymb,amsthm,amsfonts}
\usepackage[retainorgcmds]{IEEEtrantools}

\newtheorem{theorem}{Theorem}[]

\theoremstyle{definition}

\theoremstyle{remark}

\newcommand{\Zp}{{\mathbb{Z}_p}}
\newcommand{\Qp}{{\mathbb{Q}_p}}
\newcommand{\Cp}{{\mathbb{C}_p}}
\newcommand{\ZZ}{{\mathbb{Z}}}
\newcommand{\QQ}{{\mathbb{Q}}}
\newcommand{\CC}{{\mathbb{C}}}

\allowdisplaybreaks

\begin{document}

%%%%%%%%%%%%%%%%%%%%%%%%%%%%%%%%%%%%%%%%%%%%%%%%%%%%%%%%%%%%%%%%%%%%%%%%%%%%%%%%%%%%%%%%%%
%%%%%%%%%%%%%%%%%%%%%%%%%%%%%%%%%%%%%%%%%%%%%%%%%%%%%%%%%%%%%%%%%%%%%%%%%%%%%%%%%%%%%%%%%%
%%%%%%%%%%%%%%%%%%%%%%%%%%%%%%%%%%%%%%%%%%%%%%%%%%%%%%%%%%%%%%%%%%%%%%%%%%%%%%%%%%%%%%%%%%

\title[Identities for generalized twisted Bernoulli polynomials]{Identities of symmetry for generalized twisted Bernoulli polynomials twisted by unramified roots of unity}

\author{Dae San Kim}
\address{Department of Mathematics Sogang University, Seoul 121-742, South Korea}
\curraddr{} \email{dskim@sogang.ac.kr}
\thanks{}

%%%%%%%%%%%%%%%%%%%%%%%%%%%%%%%%%%%%%%%%%%%%%%%%%%%%%%%%%%%%%%%%%%%%%%%%%%%%%%%%%%%%%%%%%%
%%%%%%%%%%%%%%%%%%%%%%%%%%%%%%%%%%%%%%%%%%%%%%%%%%%%%%%%%%%%%%%%%%%%%%%%%%%%%%%%%%%%%%%%%%
%%%%%%%%%%%%%%%%%%%%%%%%%%%%%%%%%%%%%%%%%%%%%%%%%%%%%%%%%%%%%%%%%%%%%%%%%%%%%%%%%%%%%%%%%%

%\subjclass[2000]{Primary }
%    For articles to be published after 1 January 2010, you may use
%    the following version:
%\subjclass[2010]{Primary }

%%%%%%%%%%%%%%%%%%%%%%%%%%%%%%%%%%%%%%%%%%%%%%%%%%%%%%%%%%%%%%%%%%%%%%%%%%%%%%%%%%%%%%%%%%
%%%%%%%%%%%%%%%%%%%%%%%%%%%%%%%%%%%%%%%%%%%%%%%%%%%%%%%%%%%%%%%%%%%%%%%%%%%%%%%%%%%%%%%%%%
%%%%%%%%%%%%%%%%%%%%%%%%%%%%%%%%%%%%%%%%%%%%%%%%%%%%%%%%%%%%%%%%%%%%%%%%%%%%%%%%%%%%%%%%%%


%%%%%%%%%%%%%%%%%%%%%%%%%%%%%%%%%%%%%%%%%%%%%%%%%%%%%%%%%%%%%%%%%%%%%%%%%%%%%%%%%%%%%%%%%%
%%%%%%%%%%%%%%%%%%%%%%%%%%%%%%%%%%%%%%%%%%%%%%%%%%%%%%%%%%%%%%%%%%%%%%%%%%%%%%%%%%%%%%%%%%
%%%%%%%%%%%%%%%%%%%%%%%%%%%%%%%%%%%%%%%%%%%%%%%%%%%%%%%%%%%%%%%%%%%%%%%%%%%%%%%%%%%%%%%%%%

\date{}

%%%%%%%%%%%%%%%%%%%%%%%%%%%%%%%%%%%%%%%%%%%%%%%%%%%%%%%%%%%%%%%%%%%%%%%%%%%%%%%%%%%%%%%%%%
%%%%%%%%%%%%%%%%%%%%%%%%%%%%%%%%%%%%%%%%%%%%%%%%%%%%%%%%%%%%%%%%%%%%%%%%%%%%%%%%%%%%%%%%%%
%%%%%%%%%%%%%%%%%%%%%%%%%%%%%%%%%%%%%%%%%%%%%%%%%%%%%%%%%%%%%%%%%%%%%%%%%%%%%%%%%%%%%%%%%%

\dedicatory{}

%%%%%%%%%%%%%%%%%%%%%%%%%%%%%%%%%%%%%%%%%%%%%%%%%%%%%%%%%%%%%%%%%%%%%%%%%%%%%%%%%%%%%%%%%%
%%%%%%%%%%%%%%%%%%%%%%%%%%%%%%%%%%%%%%%%%%%%%%%%%%%%%%%%%%%%%%%%%%%%%%%%%%%%%%%%%%%%%%%%%%
%%%%%%%%%%%%%%%%%%%%%%%%%%%%%%%%%%%%%%%%%%%%%%%%%%%%%%%%%%%%%%%%%%%%%%%%%%%%%%%%%%%%%%%%%%

\begin{abstract}
We derive three identities of symmetry in two variables and eight in
three  variables related to generalized twisted Bernoulli
polynomials and generalized twisted power sums, both of which are
twisted by unramified roots of unity. The case of  ramified roots of
unity was treated previously. The derivations of identities are
based on the $p$-adic integral expression, with respect to a measure
introduced by Koblitz, of the generating function for the
generalized twisted Bernoulli polynomials and the quotient of
$p$-adic integrals that can be expressed as the exponential
generating function for the generalized twisted power sums.
\end{abstract}

%%%%%%%%%%%%%%%%%%%%%%%%%%%%%%%%%%%%%%%%%%%%%%%%%%%%%%%%%%%%%%%%%%%%%%%%%%%%%%%%%%%%%%%%%%
%%%%%%%%%%%%%%%%%%%%%%%%%%%%%%%%%%%%%%%%%%%%%%%%%%%%%%%%%%%%%%%%%%%%%%%%%%%%%%%%%%%%%%%%%%
%%%%%%%%%%%%%%%%%%%%%%%%%%%%%%%%%%%%%%%%%%%%%%%%%%%%%%%%%%%%%%%%%%%%%%%%%%%%%%%%%%%%%%%%%%

\maketitle

%%%%%%%%%%%%%%%%%%%%%%%%%%%%%%%%%%%%%%%%%%%%%%%%%%%%%%%%%%%%%%%%%%%%%%%%%%%%%%%%%%%%%%%%%%
%%%%%%%%%%%%%%%%%%%%%%%%%%%%%%%%%%%%%%%%%%%%%%%%%%%%%%%%%%%%%%%%%%%%%%%%%%%%%%%%%%%%%%%%%%
%%%%%%%%%%%%%%%%%%%%%%%%%%%%%%%%%%%%%%%%%%%%%%%%%%%%%%%%%%%%%%%%%%%%%%%%%%%%%%%%%%%%%%%%%%

\section{Introduction and preliminaries}\label{sec01}

Let $p$ be a fixed prime. Throughout this paper, $\Zp$, $\Qp$, $\Cp$
will respectively denote the ring of $p$-adic integers, the field of
$p$-adic rational numbers and the completion of the algebraic
closure of $\Qp$. Assume that $|\cdot|_p$ is the normalized absolute
value of $\Cp$, such that $|p|_p = \frac{1}{p}$. The group $\Gamma$
of all roots of unity of $\Cp$ is the direct product of its
subgroups $\Gamma_u$ (the subgroup of unramified roots of unity) and
$\Gamma_r$ (the subgroup of ramified roots of unity). Namely,

\begin{equation*}
\Gamma = \Gamma_{u} \cdot \Gamma_{r} ,\quad \Gamma_{u} \cap
\Gamma_{r} = \{ 1 \},
\end{equation*}

\noindent where

\begin{equation*}
\begin{split}
\Gamma_u & = \{\xi \in \CC_p | \ \xi^r = 1 \text{ for some } r \in
\ZZ_{>0} \text{ with }(r,p) = 1\},\\
\Gamma_r & = \{\xi \in \CC_p | \ \xi^{p^s} = 1 \text{ for some } s
\in \ZZ_{>0}\}.
\end{split}
\end{equation*}

Let $d$ be a fixed positive integer. Then we let

\begin{equation*}
X = X_d = \lim_{ \overleftarrow{N} } \ZZ/{dp}^N \ZZ  =
\cup_{a=0}^{dp-1}
 (a + dp \ZZ_p ),
\end{equation*}

\noindent with

\begin{equation*}
a + d p^N \ZZ_p = \{ x \in X  | \ x \equiv a \ (\bmod{\ dp^N}) \},
\end{equation*}

\noindent and let $\pi:X\rightarrow \Zp$ be the map given by the
inverse limit of the natural maps

\begin{equation*}
\ZZ/dP^N \ZZ \rightarrow \ZZ/p^N \ZZ.
\end{equation*}

If $g$ is a function on $\Zp$, we will use the same notation to
denote the function $g\circ \pi$. Let $\chi:(\ZZ/d\ZZ)^* \rightarrow
\bar{\QQ}^*$ be a (primitive) Dirichlet character of conductor $d$.
Then it will be pulled back to $X$ via the natural map $X\rightarrow
\ZZ/d\ZZ$. Here we fix, once and for all, an imbedding
$\bar{\QQ}\rightarrow \Cp$, so that $\chi$ is regarded as a map of
$X$ to $\Cp$.

Let $z\in \Cp$ be such that $z^{d p^N } \not= 1$, for all $N$. Then
we define

\begin{equation}\label{equ01}
\mu_z (a + dp^N \ZZ_p ) =  \frac{z^a}{z^{d p^N} -1}.
\end{equation}

Observe that (\ref{equ01}) is -1 times of the corresponding one in
\cite{7}. Then it is known(cf. \cite{7}) that $\mu_z$ extends to a
measure on $X$ if and only if $z\in \{ x\in \Cp \  | \ |x-1|_p \geq
1\}$. So for any such a $z$, and any continuous $\Cp$-valued
function $f$ on $X$,

{\small
\begin{equation}\label{equ02}
\begin{split}
\int_{X} {f(x) d \mu_z (x)} &= \lim_{ N \rightarrow \infty }
\sum_{a=0}
^{dp^N -1} f(a) \mu (a + dp^N \Zp )\\
 &= \lim_{ N \rightarrow \infty} \frac{1}{z^{d p^N} -1 } \sum_{a=0} ^{ dp^N -1} f(a)
 z^a.
\end{split}
\end{equation}
}

Throughout this paper, we let $\xi \not= 1$ be any fixed $r$-th root
of 1, with $(r,pd)=1$ (and hence $\xi\in\Gamma_u$), and let

\begin{equation}\label{equ03}
E = \{ t \in \Cp  \ | \  |t|_p < p^{- {\frac{1}{p-1}}} \}.
\end{equation}

Then $u_\xi$ is a measure on $X$, and, for a positive integer $w$,
$\mu_{\xi^w}$ is a measure provided that $w$ is not divisible by
$r$. For each fixed $t\in E$ the function $e^{zt}$ is analytic on
$\Zp$ and hence considered as a function on $X$. Using the
definition (\ref{equ02}), we get the $p$-adic integral expression of
the generating function for the generalized twisted Bernoulli
numbers $B_{n, \chi, \xi}$ attached to
 $\chi$ and $\xi$;

{\small
\begin{equation}\label{equ04}
t \int_{X} {\chi(z) e^{zt} d \mu_{\xi}(z)} = \frac{t}{\xi ^{d} e
^{dt} -1} \sum_{a=0}^{d-1}  \chi (a) \xi^{a} e^{at} = \sum_{n=0}
^{\infty} B_{n, \chi , \xi } \frac{t^n}{n!}  \quad (t \in E).
\end{equation}
}

\noindent So we have the following $p$-adic integral expression of
the generating function for the generalized twisted Bernoulli
polynomials $B_{n, \chi, \xi}(x)$ attached to $\chi$ and $\xi$;

{\small
\begin{equation}\label{equ05}
\begin{split}
t  \int_{X} ^{} {\chi (z)  e  ^{(x+z)t} d \mu_{\xi} (z)} &=
\frac{te^{xt}}{\xi^{d} e^{dt} -1} \sum _{a=0} ^{d-1}  \chi (a)
\xi^{a} e ^{at} \\
&= \sum _{n=0} ^{\infty} B_{n, \chi ,\xi} (x) \frac{t^{n}}{n!} \quad
(t \in E, x \in \Zp).
\end{split}
\end{equation}
}

Also, from (\ref{equ02}) we have the $p$-adic integral expression of
the generating function for the twisted Bernoulli numbers
$B_{n,\xi}$:

{\small
\begin{equation}\label{equ06}
t \int_{X }e^{zt}  d \mu_{\xi}(z) = \frac{t}{\xi e^t -1 } = \sum_{
n=0} ^{\infty} B_{n,\xi} \frac{t^n}{n!} \quad (t \in E).
\end{equation}
}

\noindent Thus we obtain the $p$-adic integral expression of the
generating function for the twisted Bernoulli polynomials $B_{n,
\xi}(x)$:

{\small
\begin{equation*}
t \int_{X } e^{(x+z)t } d \mu_{\xi}(z) = \frac{t e^{xt}}{ \xi e^t -1
} = \sum_{ n=0} ^{\infty} B_{n,\xi} (x )  \frac{ t^n}{n! } \quad (t
\in E, x \in \Zp  ).
\end{equation*}
}

Let $S_k ( n; \chi, \xi)$ denote the $k$-th generalized twisted
power sum of the first $n+1$ nonnegative integers attached to $\chi$
and $\xi$, namely
\begin{equation}\label{equ07}
S_{k}(n;\chi,\xi) = \sum_{a=0}^{n}\chi(a)\xi^{a}a^{k}=
\chi(0)\xi^00^{k}+\chi(1)\xi^11^{k}+\cdots+ \chi(n)\xi^nn^{k}.
\end{equation}

From (\ref{equ04}), (\ref{equ06}), and (\ref{equ07}), one easily
derives the following identities: for $w \in \ZZ_{> 0}$, with $w$
not divisible by $r$,

{\small
\begin{align}
\frac{\int_{X} \chi (x)e^{xt} d \mu_{\xi } (x)}{ \int_{X} e ^{d wyt}
d \mu_{\xi^{d w}} (y)} &= \frac{\xi^{dw} e^{d wt} -1}{\xi^{d}
e^{dt} -1} \sum_{a=0} ^{d-1}
\chi (a) \xi^{a} e^{at} \label{equ08}\\
&= \sum_{a=0}^{ dw -1} \chi(a)\xi^{a} e^{a t} \label{equ09}\\
&= \sum_{k=0}^{\infty} S_{k} (d w-1 ; \chi, \xi ) \frac{t^{k}}{k!}
\quad (t\in E ).\label{equ10}
\end{align}
}

In what follows, we will always assume that the $p$-adic integrals
of the various (twisted) exponential functions on $X$ are defined
for $t\in E$ (cf. (\ref{equ03})), and therefore it will not be
mentioned.

\cite{1}, \cite{2}, \cite{5}, \cite{8}, and \cite{9} are some of the
previous works on identities of symmetry in two variables involving
Bernoulli polynomials and power sums. For the brief history, one is
referred to those papers. For the first time, the idea of \cite{5}
was extended in \cite{4} to the case of three variables so as to
yield many new identities with  abundant symmetry. This added some
new identities of symmetry even to the existing ones in two
variables as well.

On the other hand, in \cite{6} the author obtained identities of
symmetry in two variables involving generalized twisted Bernoulli
polynomials and generalized twisted power sums, both of which are
twisted by ramified roots of unity (i.e., $p$-power roots of unity).
In \cite{3}, this was also extended to the case of three variables.
In these ramified cases, $p$-adic Volkenborn-type integrals are used
in both \cite{3} and \cite{6}.

In this paper, in order to treat the unramified roots of unity case
(i.e., the orders of the roots of unity are prime to $p$ and the
conductors of Dirichlet characters), we will adopt the measure
introduced by Koblitz(cf. \cite{7}) instead of Volkenborn measure,
as stated in (\ref{equ01}) and (\ref{equ02}). It seems this idea has
never been exploited before. In the end, we will be able to derive
three identities of symmetry in two variables(cf.
(\ref{equ44})-(\ref{equ46})) and eight in three variables related to
generalized twisted Bernoulli polynomials and generalized twisted
power sums (cf. (\ref{equ47})-(\ref{equ50}),
(\ref{equ53})-(\ref{equ56})).

The following is stated as Theorem \ref{thm08} and an example of the
full six symmetries in any positive integers $w_1, w_2, w_3$, with
$w_1w_2w_3$ not divisible by $r$.

{\footnotesize
\begin{align*}
& w_{1}^{n}  \sum_{k=0} ^{n} \binom{n}{k} \sum_{a=0} ^{dw_{1} -1}
\chi (a) \xi^{aw_{2} w_{3}} B_{k, \chi, \xi^{aw_{1} w_{3}}} (w_{2}
y_{1} +\frac{w_{2}}{w_{1}} a ) S_{n-k} (dw_{3} -1;
\chi, \xi^{w_{1} w_{2}} )w_{2}^{n-k} w_{3}^{k} \\
& = w_{1}^{n} \sum_{k=0} ^{n} \binom{n}{k} \sum_{a=0}^{dw_{1} -1}
\chi (a) \xi^{aw_{2} w_{3}} B_{k, \chi , \xi^{w_{1} w_{2}}} (w_{3}
y_{1} + \frac{w_{3}}{w_{1}} a ) S_{n-k} (dw_{2}
-1; \chi, \xi^{w_{1} w_{3}} )w_{3}^{n-k} w_{2}^{k} \\
& = w_{2}^{n} \sum_{k=0} ^{n} \binom{n}{k} \sum_{a=0}^{dw_{2} -1}
\chi (a) \xi^{aw_{1} w_{3}} B_{k, \chi , \xi^{w_{2} w_{3}}} (w_{1}
y_{1} + \frac{w_{1}}{w_{2}} a ) S_{n-k} (dw_{3}
-1; \chi , \xi^{w_{1} w_{2}} )w_{1}^{n-k} w_{3}^{k} \\
& = w_{2}^{n} \sum_{k=0}^{n} \binom{n}{k} \sum_{a=0} ^{dw_{2} -1}
\chi (a) \xi^{aw_{1} w_{3}} B_{k, \chi , \xi^{w_{1} w_{2}}} (w_{3}
y_{1} + \frac{w_{3}} {w_{2}} a )S_{n-k} (dw_{1}
-1; \chi ,\xi^{w_{2} w_{3}} )w_{3}^{n-k} w_{1}^{k} \\
& = w_{3}^{n} \sum_{k=0}^{n} \binom{n}{k} \sum_{a=0}^{dw_{3} -1}
\chi (a) \xi^{a w_1 w_2} B_{k, \chi , \xi^{w_2 w_3} } (w_{1} y_{1} +
\frac{w_{1}}{w_{3}} a ) S_{n-k} (dw_{2}
-1; \chi , \xi^{w_1 w_3} )w_{1}^{n-k} w_{2}^{k} \\
& = w_{3}^{n} \sum_{k=0}^{n} \binom{n}{k} \sum_{a=0}^{dw_{3} -1}
\chi (a) \xi^{aw_{1} w_{2}} B_{k,\chi , \xi^{w_{1} w_{3}}} (w_{2}
y_{1} + \frac{w_{2}}{w_{3}} a ) S_{n-k} (dw_{1} -1; \chi ,
\xi^{w_{2} w_{3}} )w_{2}^{n-k} w_{1}^{k} .
\end{align*}
}

The derivations of identities are based on the $p$-adic integral
expression of the generating function for the generalized twisted
Bernoulli polynomials in (\ref{equ05}) and the quotient of integrals
in (\ref{equ08})-(\ref{equ10}) that can be expressed as the
exponential generating function for the generalized twisted power
sums. These abundance of symmetries would not be unearthed if such
$p$-adic integral representations had not been available. We
indebted this idea to the paper \cite{5}.

%%%%%%%%%%%%%%%%%%%%%%%%%%%%%%%%%%%%%%%%%%%%%%%%%%%%%%%%%%%%%%%%%%%%%%%%%%%%%%%%%%%%%%%%%%
%%%%%%%%%%%%%%%%%%%%%%%%%%%%%%%%%%%%%%%%%%%%%%%%%%%%%%%%%%%%%%%%%%%%%%%%%%%%%%%%%%%%%%%%%%
%%%%%%%%%%%%%%%%%%%%%%%%%%%%%%%%%%%%%%%%%%%%%%%%%%%%%%%%%%%%%%%%%%%%%%%%%%%%%%%%%%%%%%%%%%

\section{Several types of quotients of $p$-adic integrals in two
variables}\label{sec02}

Here we will introduce several types of quotients of $p$-adic
integrals on $X$ or $X^2$ from which some interesting identities
follow owing to the built-in symmetries in $w_1, w_2$. In the
following, $w_1, w_2$ are positive integers with suitable
restrictions and the explicit expression of integral in
(\ref{equ11}) is obtained from the identities in (\ref{equ04}) and
(\ref{equ06}).

\vspace{3mm}

\begin{itemize}
\item[($\alpha$)] Type $\Gamma^i$ (for $i=0,1,2$)

{\footnotesize
\begin{align}
 I( \Gamma^{i} ) &= \frac{\int_{X^{2}} \chi (x_{1}
) \chi (x_{2} )t^{2-i} e^{(w_{1} x_{1} + w_{2} x_{2} + w_{1} w_{2}
(\sum_{j=1}^{2-i} y_{j} ) ) t} d \mu_{ \xi^{w_1}} (x_{1}) d
\mu_{\xi^{w_2}} (x_{2})}
{\left( \int_{X} e^{d w_1 w_2 x_3 t} d \mu_{\xi^{d w_1 w_2}} (x_3 ) \right)^{i}} \label{equ11}\\
&= \frac{{t^{2-i} e^{w_{1} w_{2} ( \sum_{j=1}^{2-i} y_{j} )t}
(\xi^{dw_{1} w_{2}} e^{dw_{1} w_{2} t }-1} )^{i} ( \sum_{a=0}^{d-1}
\chi (a) \xi^{a w_1} e^{aw_{1} t} )( \sum_{a=0}^{d-1} \chi (a)
\xi^{a w_2} e^{a w_2 t})}{( \xi^{d w_1}e^{d w_1 t} -1) (\xi^{d w_2}
e^{d w_2 t} -1) }.\label{equ12}
\end{align}
}

\end{itemize}

The above $p$-adic integral is invariant under the transposition of
$w_1,w_2$ as one can see either from the $p$-adic integral
representation in (\ref{equ11}) or from its explicit evaluation in
(\ref{equ12}).

%%%%%%%%%%%%%%%%%%%%%%%%%%%%%%%%%%%%%%%%%%%%%%%%%%%%%%%%%%%%%%%%%%%%%%%%%%%%%%%%%%%%%%%%%%
%%%%%%%%%%%%%%%%%%%%%%%%%%%%%%%%%%%%%%%%%%%%%%%%%%%%%%%%%%%%%%%%%%%%%%%%%%%%%%%%%%%%%%%%%%
%%%%%%%%%%%%%%%%%%%%%%%%%%%%%%%%%%%%%%%%%%%%%%%%%%%%%%%%%%%%%%%%%%%%%%%%%%%%%%%%%%%%%%%%%%

\section{Several types of quotients of $p$-adic integrals in three
variables}\label{sec03}

Here we will introduce several types of quotients of $p$-adic
integrals on $X$ or $X^3$ from which some interesting identities
follow due to the built-in symmetries in $w_1,w_2,w_3$. In the
following, $w_1,w_2,w_3$ are positive integers with suitable
restrictions and all of the explicit expressions of integrals in
(\ref{equ14}), (\ref{equ16}), (\ref{equ18}), and (\ref{equ20}) are
obtained from the identities in (\ref{equ04}) and (\ref{equ06}).

\vspace{3mm}

\begin{itemize}
\item[(a)] Type $\Lambda_{23}^i$ (for $i=0,1,2,3$)

{\scriptsize
\begin{align}
 I( \Lambda_{23}^{i} ) &= \frac{\substack{  \int_{X^{3}} \chi (x_{1} ) \chi (x_{2} ) \chi (x_{3}
)t^{3-i} e^{(w_{2} w_{3} x_{1} +w_{1} w_{3} x_{2} + w_{1} w_{2}
x_{3} +w_{1} w_{2} w_{3} ( \sum_{j=1}^{3-i} y_{j} ))t} \\
\hspace{4cm} \times d \mu_{\xi^{w_{2} w_{3}}} (x_{1} ) d
\mu_{\xi^{w_{1} w_{3}}} (x_{2} )d \mu_{\xi^{w_{1} w_{2}}} (x_{3} )}}
{\left( \int_{X} e^{dw_{1} w_{2} w_{3} x_{4} t} d \mu_{\xi^{dw_{1}
w_{2} w_{3}}}
(x_{4} ) \right)^{i}} \label{equ13} \\
 & = \frac{t^{3-i} e^{w_{1} w_{2} w_{3} (\sum_{j=1}^{3-i} y_{j} )t} ( \xi^{dw_{1} w_{2} w_{3}} e^{dw_{1}
w_{2} w_{3} t} -1)^{i}}{( \xi^{d w_2 w_3} e^{dw_{2} w_{3} t}
-1)(\xi^{d w_1 w_3} e^{dw_{1} w_{3} t} -1)( \xi^{d w_1 w_2}
e^{dw_{1} w_{2} t} -1)} \notag\\
 & \hspace{5mm} \times \left( \sum_{a=0}^{d-1} \chi (a) \xi ^{aw_{2} w_{3}}
e^{aw_{2} w_{3} t} \right) \left(  \sum_{a=0}^{d-1} \chi (a)
\xi^{aw_{1} w_{3}} e^{aw_{1} w_{3} t} \right) \left(
\sum_{a=0}^{d-1} \chi (a) \xi^{aw_{1} w_{2}} e^{aw_{1} w_{2} t}
\right). \label{equ14}
\end{align}
}

\noindent Here $w_2w_3,w_1w_3,w_1w_2$ are not divisible by $r$, for
$i=0$, and $w_1w_2w_3$ is not divisible by $r$, for $i=1,2,3$.

\vspace{3mm}

\item[(b)] Type $\Lambda_{13}^i$ (for $i=0,1,2,3$)

{\footnotesize
\begin{align}
I( \Lambda_{13}^{i}) & = \frac{\substack{\int_{X^{3}} \chi (x_{1} )
\chi (x_{2} ) \chi (x_{3} )t^{3-i} e^{(w_{1} x_{1} +w_{2} x_{2}
+w_{3} x_{3} +w_{1} w_{2} w_{3} ( \sum_{j=1}^{3-i} y_{j} ))t} d
\mu_{\xi^{w_{1}}} (x_{1} )  \\
\hspace{4cm} \times d \mu_{\xi^{w_{2}}} (x_{2} )d \mu_{\xi^{w_{3}}}
(x_{3} )}}{(\int_{X} e^{dw_{1} w_{2} w_{3} x_{4} t} d
\mu_{\xi^{dw_{1} w_{2} w_{3}}} (x_{4}
))^{i}} \label{equ15}\\
& = \frac{ t^{3 -i} e^{w_1 w_2 w_3 ( \sum_{ j=1}^{ 3-i} y_j ) t }
(\xi^{d w_1 w_2 w_3 }e^{dw_1 w_2 w_3 t} -1)^i }{ (\xi^{d w_1}
e^{dw_1 t} -1) (\xi^{d w_2} e^{dw_2 t} -1)(\xi^{d w_3}e^{dw_3t }
 -1)} \notag\\
 & \hspace{15mm} \times  \left( \sum_{ a=0}^{ d-1}  \chi(a ) \xi^{a w_1} e^{a w_1 t}\right)\left( \sum_{ a=0}^{ d-1} \chi(a)\xi^{a w_2} e^{a w_2 t}\right)\left( \sum_{ a=0} ^{ d-1} \chi(a) \xi^{a w_3} e^{a w_3
 t})\right).\label{equ16}
\end{align}
}

\noindent Here $w_1,w_2,w_3$ are not divisible by $r$, for $i=0$,
and $w_1w_2w_3$ is not divisible by $r$, for $i=1,2,3$.

\vspace{3mm}

\item[(c-0)] Type $\Lambda_{12}^0$

{\scriptsize
\begin{align}
I(\Lambda_{12}^{0} ) & = \int_{X^{3}} \chi (x_{1} ) \chi (x_{2} )
\chi (x_{3} )t^{3} e^{(w_{1} x_{1} +w_{2} x_{2} +w_{3} x_{3} +w_{2}
w_{3} y+w_{1} w_{3} y+w_{1} w_{2} y)t} \notag \\  & \hspace{6cm}
\times d \mu_{\xi^{w_{1}}} (x_{1} )d \mu_{\xi^{w_{2}}} (x_{2} ) d
\mu_{\xi^{w_{3}}}
(x_{3} ) \label{equ17}\\
& = \frac{ t^{3} e^{(w_{2} w_{3} +w_{1} w_{3} +w_{1} w_{2} )yt}
(\sum_{a=0}^{d-1} \chi (a) \xi^{aw_{1}} e^{aw_{1} t} )(
\sum_{a=0}^{d-1} \chi (a) \xi^{aw_{2}} e^{aw_{2} t} )(
\sum_{a=0}^{d-1} \chi (a) \xi^{aw_{3}} e^{aw_{3} t} )}{(
\xi^{dw_{1}} e^{dw_{1} t} -1)( \xi^{dw_{2}} e^{dw_{2} t} -1)(
\xi^{dw_{3}} e^{dw_{3} t} -1)}.\label{equ18}
\end{align}
}

\noindent Here $w_1,w_2,w_3$ are not divisible by $r$.

\vspace{3mm}

\item[(c-1)] Type $\Lambda_{12}^1$

{\scriptsize
\begin{align}
 I( \Lambda_{12}^{1} ) &= \frac{{ \int_{X^{3}} \chi (x_{1} ) \chi
(x_{2} ) \chi (x_{3} )e^{(w_{1} x_{1} +w_{2} x_{2} +w_{3} x_{3} )t}
d \mu_{\xi^{w_1}}(x_1 ) d \mu_{\xi^{w_2 }} (x_2 ) d \mu_{\xi^{w_3}}
(x_3 )}}{ \int_{X} e^{dw_{2} w_{3} z_{1} t} d \mu_{\xi^{dw_{2}
w_{3}}} (z_{1} ) \int_{X} e^{dw_{1} w_{3} z_{2} t} d
\mu_{\xi^{dw_{1} w_{3}}} (z_{2} ) \int_{X} e^{dw_{1} w_{2}
z_{3} t} d \mu_{\xi^{dw_{1} w_{2}}} (z_{3} )} \label{equ19}\\
& = \frac{(\xi^{d w_2 w_3 } e^{dw_{2} w_{3} t} -1)(\xi^{d w_1 w_3}
e^{dw_{1} w_{3} t} -1)(\xi^{d w_1 w_2} e^{dw_{1} w_{2} t} -1)}{(
\xi^{d w_1} e^{dw_{1} t} -1)( \xi^{d w_2} e^{dw_{2} t} -1)( \xi^{d
w_3} e^{dw_{3} t} -1)} \notag \\
 & \hspace{15mm} \times \left( \sum_{a=0}^{d-1} \chi
(a)\xi^{a w_1} e^{aw_{1} t} \right)\left( \sum_{a=0}^{d-1}  \chi (a)
\xi^{a w_2 } e^{aw_{2} t} \right)\left( \sum_{a=0}^{d-1} \chi (a
)\xi^{a w_3} e^{aw_{3} t} \right).\label{equ20}
\end{align}
}

\noindent Here $w_2w_3,w_1w_3,w_1w_2$ are not divisible by $r$.

\end{itemize}

\vspace{3mm}

All of the above $p$-adic integrals of various types are invariant
under all permutations of $w_1,w_2,w_3$, as one can see either from
$p$-adic integral representations in (\ref{equ13}), (\ref{equ15}),
(\ref{equ17}), and (\ref{equ19}) or their explicit evaluations in
(\ref{equ14}), (\ref{equ16}), (\ref{equ18}), and (\ref{equ20}).

%%%%%%%%%%%%%%%%%%%%%%%%%%%%%%%%%%%%%%%%%%%%%%%%%%%%%%%%%%%%%%%%%%%%%%%%%%%%%%%%%%%%%%%%%%
%%%%%%%%%%%%%%%%%%%%%%%%%%%%%%%%%%%%%%%%%%%%%%%%%%%%%%%%%%%%%%%%%%%%%%%%%%%%%%%%%%%%%%%%%%
%%%%%%%%%%%%%%%%%%%%%%%%%%%%%%%%%%%%%%%%%%%%%%%%%%%%%%%%%%%%%%%%%%%%%%%%%%%%%%%%%%%%%%%%%%

\section{Identities for generalized twisted Bernoulli polynomials in two
variables}\label{sec04}

All of the following results can be easily obtained from
(\ref{equ05}) and (\ref{equ08})-(\ref{equ10}).

\vspace{3mm}

\begin{itemize}
\item[($\alpha$-0)]

{\small
\begin{align}
I( \Gamma^{0}) &= \int_{X} \chi (x_{1} )te^{w_{1} (x_{1} +w_{2}
y_{1} )t} d \mu_{\xi^{w_{1}}} (x_{1} ) \int_{X} \chi (x_{2}
)te^{w_{2} (x_{2} +w_{1} y_{2} )t} d \mu_{\xi^{w_{2}}} (x_{2} ) \notag\\
& = \left( \sum_{k=0}^{\infty} B_{k, \chi , \xi^{w_{1}}} (w_{2}
y_{1} )
\frac{(w_{1} t)^{k}}{k!} \right)\left( \sum_{l=0}^{\infty} B_{l, \chi, \xi^{w_2}} (w_1 y_2 ) \frac{ (w_2 t)^l}{ l!} \right) \notag\\
& = \sum_{ n=0}^{\infty} \left( \sum_{ k=0}^{ n} \binom{n}{k}B_{k,
\chi, \xi^{w_1}} (w_2 y_1 ) B_{n-k, \chi, \xi^{w_2 }}(w_1 y_2 )
w_1^{k} w_2^{n-k} \right)  \frac{t^n}{n!}.\label{equ21}
\end{align}
}

\item[($\alpha$-1)] Here we write $I( \Gamma^{1} )$ in two different ways:

\vspace{5mm}
\item[(1)]

{\small
\begin{align}
I( \Gamma^{1} ) &= \int_{X} \chi (x_{1} )te^{w_{1} (x_{1} +w_{2}
y_{1} )t} d \mu_{\xi^{w_{1}}} (x_{1} ) \times \frac{ \int_{X} \chi
(x_{2} )e^{w_{2} x_{2} t} d \mu_{\xi^{w_{2}}}(x_2 )}{ \int_{ X} e^{
d w_1 w_2 x_3 t}d \mu_{\xi^{d w_1 w_2}}(x_3
 )} \label{equ22}\\
& = \left( \sum_{ k=0} ^{ \infty} B_{k, \chi, \xi^{w_1}} (w_2 y_1 )
\frac{ (w_1 t)^k}{ k!} \right) \left( \sum_{ l=0}^{ \infty} S_l
(dw_1 -1; \chi,
\xi^{w_2}) \frac{ (w_2 t)^l}{ l!} \right) \notag\\
 & = \sum_{ n=0}^{ \infty} \left(\sum_{ k=0}^{n}\binom{n}{k}B_{k,\chi,\xi^{w_1}} (w_2 y_1 ) S_{n-k} (d w_1
-1;\chi, \xi^{w_2} ) w_1^k w_2^{n-k} \right) \frac{ t^n}{
n!}.\label{equ23}
\end{align}
}

\vspace{5mm}
\item[(2)] Invoking (\ref{equ09}), (\ref{equ22}) can also be written as

{\small
\begin{align}
I( \Gamma^{1} ) &= \sum_{a=0}^{ dw_1 -1} \chi (a ) \xi^{a w_2}
\int_{ X} \chi (x_1 ) t e^{w_1 ( x_1 +w_2 y_1 + \frac{ w_2}{ w_1}a ) t}d \mu_{\xi^{w_1}} (x_1 ) \notag\\
& = \sum_{ a=0}^{ dw_1 -1} \chi (a ) \xi^{a w_2} \left( \sum_{
n=0}^{ \infty} B_{n, \chi, \xi^{w_1}} (w_2 y_1 + \frac{ w_2}{w_1 }
a)
\frac{(w_1 t)^n}{n! } \right) \notag\\
& = \sum_{ n=0}^{\infty} \left(w_1^{n} \sum_{ a=0}^{ dw_1 -1} \chi
(a) \xi^{a w_2} B_{n, \chi, \xi^{w_1}}(w_2 y_1 + \frac{ w_2}{w_1}a)
\right) \frac{ t^n}{ n!}.\label{equ24}
\end{align}
}

\vspace{5mm}
\item[($\alpha$-2)]

{\small
\begin{align*}
I( \Gamma^{2} ) &= \frac{ \int_{X} \chi (x_{1} ) e^{w_{1} x_{1} t} d
\mu_{\xi^{w_{1}}} (x_{1} )}{ \int_{X} e^{dw_{1} w_{2} x_{3} t} d
\mu_{\xi^{dw_{1} w_{2}}} (x_{3} )} \times  \frac{ \int_{X} \chi
(x_{2} ) e^{w_{2} x_{2} t} d \mu_{ \xi^{w_{2}}}(x_2 )}{ \int_{
X}  e^{d w_1 w_2 x_3 t} d \mu_{\xi^{ d w_1 w_2 }}( x_3 ) }\\
& = \left( \sum_{k=0}^{\infty} S_{k} (dw_{2} -1; \chi , \xi^{w_{1}}
) \frac{(w_1 t )^k}{k!} \right) \left( \sum_{l=0}^{ \infty} S_l (d
w_1-1; \chi, \xi^{ w_2 } ) \frac{ (w_2 t)^l}{l! } \right)\\
& = \sum_{n=0}^{\infty} \left( \sum_{k=0}^{n} \binom{n}{k} S_{k}
(dw_{2} -1; \chi , \xi^{w_{1}} )S_{n-k} (dw_{1} -1; \chi ,
\xi^{w_{2}} ) w_{1}^{k} w_{2}^{n-k} \right) \frac{t^{n}}{n!}.
\end{align*}
}

\end{itemize}

%%%%%%%%%%%%%%%%%%%%%%%%%%%%%%%%%%%%%%%%%%%%%%%%%%%%%%%%%%%%%%%%%%%%%%%%%%%%%%%%%%%%%%%%%%
%%%%%%%%%%%%%%%%%%%%%%%%%%%%%%%%%%%%%%%%%%%%%%%%%%%%%%%%%%%%%%%%%%%%%%%%%%%%%%%%%%%%%%%%%%
%%%%%%%%%%%%%%%%%%%%%%%%%%%%%%%%%%%%%%%%%%%%%%%%%%%%%%%%%%%%%%%%%%%%%%%%%%%%%%%%%%%%%%%%%%

\section{Identities for generalized twisted Bernoulli polynomials in three
variables}\label{sec05}

All of the following results can be easily obtained from
(\ref{equ05}) and (\ref{equ08})-(\ref{equ10}). First, let's consider
Type $\Lambda_{23}^i$, for each $i=0,1,2,3$.

\begin{itemize}

\vspace{5mm}
\item[(a-0)]

{\scriptsize
\begin{align}
I( \Lambda_{23}^0 ) &= \int_{X}  {\chi (x_1 )te^{w_2 w_3 (x_1+w_1 y_1 )t }d
\mu_{\xi^{w_2 w_3}}}( x_1 ) \int_{ X}  {\chi (x_2)te^{w_1 w_3
(x_2+w_2 y_2 ) t}d \mu_{\xi^{w_1 w_3}}} ( x_2 ) \notag \\
& \hspace{5cm} \times \int_{ X}{\chi(x_3 )t e^{w_1 w_2 (x_3 +w_3 y_3 ) t }d \mu_{\xi^{w_1 w_2}}} (x_3 )\notag\\
& = \left(\sum_{k=0}^{\infty}  \frac{B_{k,\chi ,\xi^{w_{2} w_{3}}}
(w_{1} y_{1} )}{k!} (w_{2} w_{3} t)^{k} \right)\left(
\sum_{l=0}^{\infty} \frac{B_{l, \chi , \xi^{w_{1} w_{3}}} (w_{2}
y_{2} )}{l!} (w_{1} w_{3} t)^{l} \right) \notag \\
 & \hspace{5cm}
\times \left( \sum_{m=0} ^{\infty}
\frac{B_{m, \chi ,\xi^{w_{1} w_{2}}} (w_{3} y_{3})}{m!} (w_{1} w_{2} t)^{m} \right) \notag \\
 & = \sum_{n=0}^{\infty} ( \sum_{k+l+m=n} \binom{n}{k,l,m} B_{k,\chi ,\xi^{w_{2} w_{3}}} (w_{1} y_{1} )B_{l, \chi ,
\xi^{w_{1} w_{3}}} (w_{2} y_{2} ) B_{m, \chi , \xi^{w_{1} w_{2}}}
(w_{3} y_{3} ) \notag \\
 & \hspace{5cm} \times  w_{1}^{l+m}
w_{2}^{k+m} w_{3}^{k+l}) \frac{ t^n}{n! },\label{equ25}
\end{align}
}

\noindent where the inner sum is over all nonnegative integers
$k,l,m$, with $k+l+m=n$ and

{\small
\begin{align*}
\binom{n}{k,l,m} = \frac{ n!}{ k!\  l!\  m!}.
\end{align*}
}

\item[(a-1)] Here we write $I( \Lambda_{23}^{1} )$ in two different ways:

\vspace{5mm}
\item[(1)]

{\scriptsize
\begin{align}
I( \Lambda_{23}^{1} ) &= \int_{X} {\chi (x_{1}
)te^{w_{2} w_{3} (x_{1} +w_{1} y_{1} )t} d \mu_{\xi^{w_{2} w_{3}}}}
(x_{1} ) \int_{X} { \chi (x_{2} )te^{w_{1} w_{3} (x_{2} +w_{2} y_{2}
)t} d \mu_{\xi^{w_{1} w_{3}}}} (x_{2} ) \notag \\
 & \hspace{5cm} \times
\frac{ \int_{ X} \chi
(x_3 )e^{w_1 w_2 x_3 t} d \mu_{\xi^{w_1 w_2 }} ( x_3 ) }{\int_{X} {e^{dw_1 w_2 w_3 x_4t}d  \mu_{\xi^{d w_1 w_2 w_3 }} (x_4  ) } } \label{equ26}\\
& = \left( \sum_{k=0}^{\infty} B_{k, \chi , \xi^{w_{2} w_{3}}}
(w_{1} y_{1} ) \frac{(w_{2} w_{3} t)^{k}}{k!} \right)\left(
\sum_{l=0}^{\infty} B_{l, \chi , \xi^{w_{1} w_{3}}} (w_{2} y_{2} )
\frac{(w_{1} w_{3} t)^{l}}{l!} \right) \notag \\ & \hspace{5cm}
\times \left(\sum_{ m=0}^{\infty} S_{m} (d w_{3} -1;
\chi, \xi^{w_1w_2} ) \frac{(w_{1} w_{2} t)^{m}}{m!} \right)\notag\\
& = \sum_{n=0}^{\infty} ( \sum_{k+l+m=n} \binom{n}{k,l,m} B_{k, \chi
, \xi^{w_2 w_3} } (w_{1} y_{1} )B_{l, \chi , \xi^{w_1 w_3} } (w_{2}
y_{2} )  S_{m} (d w_{3} -1; \chi , \xi^{w_1w_2})\notag \\ &
\hspace{5cm} \times w_{1}^{l+m} w_{2}^{k+m} w_{3}^{k+l} ) \frac{
t^n}{n! }.\label{equ27}
\end{align}
}

\vspace{5mm}
\item[(2)] Invoking (\ref{equ09}), (\ref{equ20}) can also be written as

{\scriptsize
\begin{align}
 I( \Lambda_{23}^{1} ) &= \sum_{a=0}^{d w_{3} -1}
\chi (a) \xi^{a w_1 w_2} \int_{X} \chi (x_{1} )te^{w_{2} w_{3}
(x_{1} +w_{1} y_{1} )t} d \mu_{\xi^{ w_2 w_3 }} (x_{1}) \notag \\
& \hspace{5cm} \times \int_{X} {\chi (x_{2} )te^{w_{1} w_{3} (x_{2}
+w_{2} y_{2} + \frac{w_{2}}{w_{3}}
a)t} d \mu_{\xi^{w_1 w_3} } (x_{2} )} \notag\\
 & = \sum_{a=0}^{dw_{3} -1} \chi (a) \xi^{aw_{1} w_{2}} \left(\sum_{k=0}^{\infty } B_{k, \chi , \xi^{w_{2} w_{3}}} (w_{1} y_{1} ) \frac{(w_{2} w_{3} t)^{k}}{k!} \right) \notag \\
 & \hspace{5cm} \times \left( \sum_{l=0}^{\infty } B_{l, \chi , \xi^{w_{1} w_{3}}} (w_{2} y_{2}
+ \frac{w_{2}}{w_{3}} a) \frac{(w_{1} w_{3} t)^{l}}{l!} \right)\notag \\
 & = \sum_{n=0}^{\infty} (w_{3}^{n} \sum_{k=0}^{n} \binom{n}{k} B_{k, \chi , \xi^{w_{2} w_{3}}} (w_{1} y_{1} )
\sum_{a=0}^{dw_{3} -1} \chi (a) \xi^{aw_{1} w_{2}} \notag \\
& \hspace{5cm} \times B_{n-k, \chi , \xi^{w_{1} w_{3}}} (w_{2} y_{2}
+ \frac{w_{2}}{w_{3}} a)w_{1}^{n-k} w_{2}^{k} )
\frac{t^{n}}{n!}.\label{equ28}
\end{align}
}

\vspace{5mm}
\item[(a-2)] Here we write $I( \Lambda_{23}^{2})$ in three different ways:

\vspace{5mm}
\item[(1)]

{\tiny
\begin{align}
I( \Lambda_{23}^{2} ) &= \int_{X} {\chi (x_{1} )te^{w_{2} w_{3} (x_{1} + w_{1} y_{1} )t}
d \mu_{ \xi^{w_{2} w_{3}}}} (x_{1} ) \notag \\
 & \qquad  \times \frac{
\int_{X} \chi (x_{2} )e^{w_{1} w_{3} x_{2} t} d \mu_{\xi^{w_{1}
w_{3}}} (x_{2} )}{ \int_{X} e^{dw_{1} w_{2} w_{3} x_{4} t} d
\mu_{\xi^{dw_{1} w_{2} w_{3}}} (x_{4} )} \times \frac{ \int_{X}
{\chi (x_{3} ) e^{w_{1} w_{2} x_{3} t}} d \mu_{\xi^{ w_1 w_2} }
(x_{3} )}{ \int_{X} {e^{dw_{1} w_{2}
w_{3} x_{4} t} d \mu_{\xi^{dw_{1} w_{2} w_{3}}} (x_{4} )}} \label{equ29}\\
& = \left( \sum_{k=0}^{\infty} B_{k, \chi , \xi^{w_2 w_3} } (w_{1}
y_{1} ) \frac{(w_{2} w_{3} t)^{k}}{k!} \right)\left(
\sum_{l=0}^{\infty} S_{l} (dw_{2} -1; \chi , \xi^{w_1 w_3} )
\frac{(w_{1} w_{3} t)^{l}}{l!} \right) \notag \\
 & \hspace{5cm} \times \left(
\sum_{ m=0}^{\infty } S_m (dw_3
-1; \chi, \xi^{w_1 w_2}) \frac{ (w_1 w_2 t)^m}{m! } \right)\notag \\
& = \sum_{n=0}^{\infty} ( \sum_{k+l+m=n} \binom{n}{k,l,m} B_{k, \chi
, \xi^{w_{2} w_{3}}} (w_{1} y_{1} )S_{l} (dw_{2} -1; \chi
,\xi^{w_{1} w_{3}} )S_{m} (dw_{3} -1; \chi , \xi^{w_{1} w_{2}} )
\notag \\
& \hspace{5cm} \times
 w_{1}^{l+m} w_{2}^{k+m} w_{3}^{k+l} )
\frac{t^{n}}{n!}.\label{equ30}
\end{align}
}

\vspace{5mm}
\item[(2)] Invoking (\ref{equ09}), (\ref{equ29}) can also be written as

{\scriptsize
\begin{align}
I( \Lambda_{23}^{2} ) &= \sum_{a=0}^{dw_{2} -1} \chi (a) \xi^{aw_{1}
w_{3}} \int_{X} \chi (x_{1} )te^{w_{2} w_{3} (x_{1} +w_{1} y_{1} +
\frac{w_{1}}{w_{2}} a)t} d \mu_{\xi^{ w_2 w_3} } (x_{1} ) \notag \\
& \hspace{5cm} \times \frac{ \int_{X} {\chi (x _{3} )e^{w_{1} w_{2}
x_{3} t}} d \mu_{\xi^{w_{1} w_{2}}} (x_{3} )}{ \int_{X} {e^{dw_{1}
w_{2}
w_{3} x_{4} t} d \mu_{\xi^{dw_{1} w_{2} w_{3}}} (x_{4} )}}\label{equ31}\\
& = \sum_{a=0}^{dw_{2} -1} \chi (a ) \xi^{a w_1 w_3} \left(
\sum_{k=0}^{\infty} B_{k, \chi, \xi^{w_2 w_3}} (w_{1} y_{1} +
\frac{w_{1}}{w_{2}} a) \frac{(w_{2} w_{3} t)^{k}}{k!} \right)
\notag\\
& \hspace{5cm} \times \left( \sum_{l=0}^{\infty} S_{l} (d w_{3} -1;
\chi, \xi^{w_1 w_2 })
\frac{(w_{1} w_{2} t)^{l}}{l!} \right)\notag\\
& = \sum_{n=0}^{\infty} (w_{2}^{n} \sum_{k=0}^{n} \binom{n}{k}
\sum_{a=0}^{dw_{2} -1} \chi (a) \xi^{aw_{1} w_{3}} B_{k, \chi
,\xi^{w_{2} w_{3}}} (w_{1} y_{1} +
\frac{w_{1}}{w_{2}} a)\notag\\
& \hspace{5cm} \times S_{n-k} (dw_{3} -1; \chi , \xi^{w_{1} w_{2}}
)w_{1}^{n-k} w_{3}^{k} ) \frac{t^{n}}{n!}.\label{equ32}
\end{align}
}

\item[(3)] Invoking (\ref{equ09}) once again, (\ref{equ31}) can be written
as

{\scriptsize
\begin{align}
I( \Lambda_{23}^{2} ) &= \sum_{a=0}^{dw_{2} -1} \chi
(a) \xi^{a w_1 w_3} \sum_{b=0}^{dw_{3} -1} \chi (b) \xi^{bw_1 w_2}
\int_{X} {\chi
(x _{1} )te^{w_{2} w_{3} (x_{1} +w_{1} y_{1} + \frac{w_{1}}{w_{2}} a+ \frac{w_{1}}{w_{3}} b)t} d \mu_{\xi^{ w_2 w_3 } } (x_{1} )} \notag\\
& = \sum_{ a=0}^{ dw_2 -1 } \chi(a ) \xi^{a w_1 w_3} \sum_{
b=0}^{dw_3 -1 } \chi (b ) \xi^{b w_1 w_2 } \left(\sum_{ n=0}^{
\infty} B_{n, \chi, \xi^{ w_2 w_3}} (w_1 y_1+ \frac{ w_1}{ w_2}a+
\frac{w_1}{ w_3}b ) \frac{(w_2 w_3 t)^n}{ n!}\right) \notag \\
& = \sum_{n=0}^{\infty } \left((w_{2} w_{3} )^{n}
\sum_{a=0}^{dw_{2} -1}
 \sum_{b=0}^{dw_{3} -1} \chi (ab) \xi^{w_{1} (aw _{3} +bw_{2} )} B_{n,\chi , \xi^{w_{2} w_{3}}} (w_{1} y_{1}
+ \frac{w_{1}}{w_{2}} a+ \frac{w_{1}}{w_{3}} b)\right)
\frac{t^{n}}{n!}.\label{equ33}
\end{align}
}

\vspace{5mm}
\item[(a-3)]

{\scriptsize
\begin{align}
I(\Lambda_{23}^{3}) &= \frac{\int_{X}\chi(x_{1})
e^{w_{2}w_{3}x_{1}t}d\mu_{\xi^{w_{2}w_{3}}}(x_{1})}{\int_{X}
e^{dw_{1}w_{2}w_{3}x_{4}t}d\mu_{\xi^{dw_{1}w_{2}w_{3}}} (x_{4})}
\times \frac{\int_{X} \chi(x_{2})e^{w_{1}w_{3}x_{2}t}
d\mu_{\xi^{w_{1}w_{3}}}(x_{2})}{\int_{X} e^{dw_{1}w_{2}w_{3}x_{4}t}
d \mu_{\xi^{dw_{1}w_{2}w_{3}}}(x_{4})} \notag\\
& \hspace{5cm} \times\frac{\int_{X} \chi(x_{3})e^{w_{1}w_{2}x_{3}t}
d\mu_{\xi^{w_{1}w_{2}}}(x_{3})}{\int_{X} e^{dw_1w_2
w_3x_4t} d\mu_{\xi^{dw_1w_2w_3}}(x_4)}\notag\\
& =\left(\sum_{k=0}^{\infty}S_{k}(dw_{1}-1;\chi,\xi^{w_{2}w_{3}})
\frac{(w_{2}w_{3}t)^{k}}{k!}\right)\left(\sum_{l=0}^{\infty}S_{l}(dw_{2}-1;
\chi,\xi^{w_1w_3}) \frac{(w_{1}w_{3}t)^{l}}{l!}\right) \notag \\
& \hspace{5cm} \times
\left(\sum_{m=0}^{\infty}  S_m(dw_3-1;\chi,\xi^{w_1w_2}) \frac{(w_1w_2t)^m}{m!}\right)\notag\\
&= \sum_{n=0}^{\infty} (\sum_{k+l+m=n}\binom{n}{k,l,m} S_k
(dw_1-1;\chi,\xi^{w_2w_3})S_l(dw_2-1;\chi,\xi^{w_1 w_3}) \notag\\
& \hspace{5cm} \times S_m(dw_3-1;\chi,\xi^{w_1w_2}) w_1^{l+m}
w_2^{k+m}w_3^{k+l})\frac{t^n}{n!}.\label{equ34}
\end{align}
}

\item[(b)] For Type $\Lambda_{13}^i$ ($i=0,1,2,3$), we may consider the
analogous things to the ones in (a-0), (a-1), (a-2), and (a-3).
However, each of those can be obtained from the corresponding ones
in (a-0), (a-1), (a-2), and (a-3).  Indeed, if we substitute $w_2
w_3,w_1 w_3,w_1 w_2$ respectively for $w_1,w_2,w_3$ in $\frac{I(
\Lambda_{23}^i )}{t^{3-i}}$ (cf. (\ref{equ13})), this amounts to
replacing $t$ by $w_1 w_2 w_3 t$ and $\xi$ by $\xi^{w_1 w_2 w_3}$ in
$\frac{I( \Lambda_{13}^i )}{t^{3-i}}$ (cf. (\ref{equ15})). So, upon
replacing $w_1,w_2,w_3$ respectively by $w_2 w_3,w_1 w_3,w_1 w_2$ ,
dividing by $(w_1 w_2 w_3 )^n$, and replacing $\xi^{w_1 w_2 w_3}$ by
$\xi$, in each of the expressions of (\ref{equ25}), (\ref{equ27}),
(\ref{equ28}), (\ref{equ30}), (\ref{equ32})-(\ref{equ34}), we will
get the corresponding symmetric identities for Type $\Lambda_{13}^i$
($i=0,1,2,3$).

\vspace{5mm}
\item[(c-0)]

{\scriptsize
\begin{align}
I(\Lambda_{12}^0) &= \int_{X} \chi(x_{1})te^{w_{1} (x_{1}+w_{2}y)t}d
\mu_{\xi^{w_{1}}}(x_{1}) \int_{X} \chi(x_{2})te^{w_{2}(x_{2}+w_{3}
y)t} d\mu_{\xi^{w_{2}}}(x_{2}) \notag\\
& \hspace{5cm} \times \int_{X}
\chi(x_{3})te^{w_{3}(x_{3}+w_{1}y)t} d\mu_{\xi^{w_3}}(x_{3}) \notag\\
&=
\left(\sum_{k=0}^{\infty}\frac{B_{k,\chi,\xi^{w_1}}(w_2y)}{k!}(w_1t)^k \right)
\left(\sum_{l=0}^{\infty}\frac{B_{l,\chi,\xi^{w_2}}(w_3
y)}{l!}(w_2t)^l\right)\left(\sum_{m=0}^{\infty}\frac{B_{m,\chi,\xi^{w_3}} (w_1y)}{m!}(w_3t)^m \right) \notag\\
&= \sum_{n=0}^{\infty}\left(\sum_{k+l+m=n}\binom{n}{k,l,m}
B_{k,\chi,\xi^{w_1}}(w_2y)B_{l,\chi,\xi^{w_2}}(w_3y)
B_{m,\chi,\xi^{w_3}}(w_1y) w_1^k w_2^l w_3^m  \right) \frac{t^n}{n!
}.\label{equ35}
\end{align}
}

\item[(c-1)]

{\scriptsize
\begin{align}
I(\Lambda_{12}^1)&=\frac{{\int_{X}\chi(x_1)e^{w_1x_1t}
d\mu_{\xi^{w_1}}(x_1)}}{\int_{X}e^{dw_1w_2z_3t}d
\mu_{\xi^{dw_1w_2}}(z_3)}\times\frac{\int_{X}\chi(x_2
)e^{w_2x_2t}d\mu_{\xi^{w_2}}(x_2)}{\int_{X}e^{dw_2w_3
z_1t}d\mu_{\xi^{dw_2w_3}}(z_1)}\times \frac{\int_{X}\chi
(x_3)e^{w_3x_3t}d\mu_{\xi^{w_3}}(x_3)}{\int_{X}
e^{dw_3w_1z_2t}d\mu_{\xi^{dw_3w_1}}(z_2)} \notag \\
&=\left(\sum_{k=0}^{\infty}S_k(dw_2-1;\chi,\xi^{w_1})\frac{(w_1
t)^k}{k!}\right)\left(\sum_{l=0}^{\infty}S_l(dw_3-1;\chi, \xi^{w_2})
\frac{(w_2t)^l}{l!}\right) \notag
\\ & \hspace{5cm} \times \left(\sum_{m=0}^{\infty}
S_m(dw_1-1;\chi,\xi^{w_3})\frac{(w_3t)^m}{m!}\right) \notag \\
&=\sum_{n=0}^{\infty}(\sum_{k+l+m=n} \binom{n}{k,l,m}
S_{k}(dw_2-1;\chi,\xi^{w_1})S_l(dw_3-1;\chi, \xi^{w_2}) \notag
\\
& \hspace{5cm} \times
S_m(dw_1-1;\chi,\xi^{w_3})w_1^{k}w_2^{l}w_3^{m})\frac{t^n}{n!}.\label{equ36}
\end{align}
}

\end{itemize}

%%%%%%%%%%%%%%%%%%%%%%%%%%%%%%%%%%%%%%%%%%%%%%%%%%%%%%%%%%%%%%%%%%%%%%%%%%%%%%%%%%%%%%%%%%
%%%%%%%%%%%%%%%%%%%%%%%%%%%%%%%%%%%%%%%%%%%%%%%%%%%%%%%%%%%%%%%%%%%%%%%%%%%%%%%%%%%%%%%%%%
%%%%%%%%%%%%%%%%%%%%%%%%%%%%%%%%%%%%%%%%%%%%%%%%%%%%%%%%%%%%%%%%%%%%%%%%%%%%%%%%%%%%%%%%%%

\section{Main theorems}\label{sec06}

As we noted earlier, the various types of quotients of $p$-adic
integrals in Section \ref{sec02} and Section \ref{sec03} are
respectively invariant under the transposition of $w_1,w_2$ and any
permutation of $w_1,w_2,w_3$. So the corresponding expressions in
Section \ref{sec04} and Section \ref{sec05} are respectively also
invariant under the transposition of $w_1,w_2$ and any permutation
of $w_1,w_2,w_3$. Thus our results about identities of symmetry will
be immediate consequences of this observation.

However, not all permutations of an expression in Section
\ref{sec05} yield distinct ones. In fact, as these expressions are
obtained by permuting $w_1,w_2,w_3$ in a single one labeled by them,
they can be viewed as a group in a natural manner and hence it is
isomorphic to a quotient of $S_3$. In particular, the number of
possible distinct expressions are 1,2,3, or 6. (a-0), (a-1(1)),
(a-1(2)), and (a-2(2)) give the full six identities of symmetry,
(a-2(1)) and (a-2(3)) yield three identities of symmetry, and (c-0)
and (c-1) give two identities of symmetry, while the expression in
(a-3) yields no identities of symmetry. Similarly, ($\alpha$-0),
($\alpha$-1(1)), and ($\alpha$-1(2)) give two identities of symmetry
but ($\alpha$-2) yields no identity of symmetry.

Here we will just consider the cases of Theorems \ref{thm07} and
\ref{thm11}, leaving the others as easy exercises for the reader. As
for the case of Theorem \ref{thm07}, in addition to
(\ref{equ50})-(\ref{equ52}), we get the following three ones:

{\scriptsize
\begin{align}
&\sum_{k+l+m=n}\binom{n}{k,l,m} B_{k,\chi,\xi^{w_{2}w_{3}}}
(w_{1}y_{1})S_{l}(dw_{3}-1;\chi
,\xi^{w_{1}w_{2}})S_{m}(dw_{2}-1;\chi,\xi^{w_{1}w_{3}})
w_{1}^{l+m}w_{3}^{k+m}w_{2}^{k+l},\label{equ37}\\
&\sum_{k+l+m=n}\binom{n}{k,l,m}B_{k,\chi,\xi^{w_{1}w_{3}}}
(w_{2}y_{1})S_{l}(dw_{1}-1;\chi,
\xi^{w_{2}w_{3}})S_{m}(dw_{3}-1;\chi,\xi^{w_{1}
w_{2}})w_{2}^{l+m}w_{1}^{k+m}w_{3}^{k+l},\label{equ38}\\
&\sum_{k+l+m=n}\binom{n}{k,l,m}B_{k,\chi,\xi^{w_{1}w_{2}}}
(w_{3}y_{1})S_{l}(dw_{2}-1;\chi,
x^{w_{1}w_{3}})S_{m}(dw_{1}-1;\chi,\xi^{w_{2}
w_{3}})w_{3}^{l+m}w_{2}^{k+m}w_{1}^{k+l}.\label{equ39}
\end{align}
}

But, by interchanging $l$ and $m$, we see that (\ref{equ37}),
(\ref{equ38}), and (\ref{equ39}) are respectively equal to
(\ref{equ50}), (\ref{equ51}), and (\ref{equ52}). As to Theorem
\ref{thm11}, in addition to (\ref{equ56}) and (\ref{equ57}), we
have:

{\footnotesize
\begin{align}
&\sum_{k+l+m=n}\binom{n}{k,l,m}
S_k(dw_2-1;\chi,\xi^{w_1})S_l(dw_3-1;\chi,\xi^{w_2})
S_m(dw_1-1;\chi,\xi^{w_3})w_1^{k}w_2^{l}w_3^{m},\label{equ40}\\
&\sum_{k+l+m=n}\binom{n}{k,l,m}S_{k}
(dw_{3}-1;\chi,\xi^{w_{2}})S_{l}(dw_{1}-1;\chi,
\xi^{w_{3}})S_{m}(dw_{2}-1;\chi,\xi^{w_{1}})w_2^{k}w_3^{l}w_1^{m},\label{equ41}\\
&\sum_{k+l+m=n}\binom{n}{k,l,m} S_{k}(dw_{3}-1; \chi
,\xi^{w_1})S_{l}(dw_{2}-1;\chi,\xi^{w_3}) S_{m}
(dw_{1}-1;\chi,\xi^{w_2})w_1^{k}w_3^{l}w_2^{m},\label{equ42}\\
&\sum_{k+l+m=n}\binom{n}{k,l,m}S_{k}(dw_{2}-1;\chi
,\xi^{w_3})S_{l}(dw_{1}-1;\chi,\xi^{w_2})S_{m}
(dw_{3}-1;\chi,\xi^{w_1})w_3^{k}w_2^{l}w_1^{m}.\label{equ43}
\end{align}
}

However, (\ref{equ40}) and (\ref{equ41}) are equal to (\ref{equ56}),
as we can see by applying the permutations $k \rightarrow l,
l\rightarrow m, m\rightarrow k$ for (\ref{equ40}) and $k \rightarrow
m, l\rightarrow k, m\rightarrow l$ for (\ref{equ41}). Similarly, we
see that (\ref{equ42}) and (\ref{equ43}) are equal to (\ref{equ57}),
by applying permutations $k\rightarrow l, l\rightarrow m, ,m
\rightarrow k$ for (\ref{equ42}) and $k\rightarrow m, l\rightarrow
k, m\rightarrow l$ for (\ref{equ43}).

\begin{theorem}\label{thm01}
Let $w_1,w_2$ be any positive integers, such that $r$ does not
divide $w_1,w_2$. Then we have: {\small
\begin{equation}\label{equ44}
\begin{split}
& \sum_{k=0}^{n}\binom{n}{k}B_{k,\chi,\xi^{w_1}}(w_2y_1)B_{n-k,\chi,
\xi^{w_2}}(w_1y_2)w_1^{k} w_2^{n-k} \\ & \quad =\sum_{k=0}^{n}\binom{n}{k}
B_{k,\chi,\xi^{w_2}}(w_1y_1)B_{n-k,\chi,\xi^{w_1}}(w_2y_2)w_2^{k}w_1^{n-k}.
\end{split}
\end{equation}
}
\end{theorem}

\begin{theorem}\label{thm02}
Let $w_1,w_2$ be any positive integers, such that $r$ does not
divide $w_1w_2$. Then we have: {\small
\begin{equation}\label{equ45}
\begin{split}
&
\sum_{k=0}^{n}\binom{n}{k}B_{k,\chi,\xi^{w_1}}(w_2y_1)S_{n-k}(dw_1-1;\chi,
\xi^{w_2})w_1^kw_2^{n-k} \\ & \quad = \sum_{k=0}^{n}\binom{n}{k}
B_{k,\chi,\xi^{w_2}}(w_1y_1)
S_{n-k}(dw_2-1;\chi,\xi^{w_1})w_2^kw_1^{n-k}.
\end{split}
\end{equation}
}
\end{theorem}

\begin{theorem}\label{thm03}
Let $w_1,w_2$ be any positive integers, such that $r$ does not
divide $w_1w_2$. Then we have: {\small
\begin{equation}\label{equ46}
\begin{split}
& w_1^{n}\sum_{a=0}^{dw_1-1}\chi(a)\xi^{aw_2}
B_{n,\chi,\xi^{w_1}}(w_2y_1+\frac{w_2}{w_1}a) \\ & \quad =w_2^{n}
\sum_{a=0}^{dw_2-1}\chi(a)\xi^{aw_1}B_{n,\chi,\xi^{w_2}}(w_1
y_1+\frac{w_1}{w_2}a).
\end{split}
\end{equation}
}
\end{theorem}

\begin{theorem}\label{thm04}
Let $w_1,w_2,w_3$ be any positive integers, such that $r$ does not
divide $w_2w_3,w_1w_3,w_1w_2$. Then we have: {\footnotesize
\begin{equation}\label{equ47}
\begin{split}
&\sum_{k+l+m=n} \binom{n}{k,l,m}B_{k,\chi,
\xi^{w_{2}w_{3}}}(w_{1}y_{1})B_{l,\chi,\xi^{w_{1} w_{3}}} (w_{2}y_{2})B_{m,\chi, \xi^{w_{1}w_{2}}} (w_{3}y_{3})w_{1}^{l+m}w_{2}^{k+m} w_{3}^{k+l}\\
&=\sum_{k+l+m=n} \binom{n}{k,l,m} B_{k,\chi,
\xi^{w_{2}w_{3}}}(w_{1}y_{1})B_{l,\chi,
\xi^{w_{1}w_{2}}}(w_{3}y_{2})
B_{m,\chi,\xi^{w_{1}w_{3}}} (w_{2}y_{3})w_{1}^{l+m} w_{3}^{k+m} w_{2}^{k+l}\\
&= \sum_{k+l+m=n} \binom{n}{k,l,m} B_{k,\chi,
\xi^{w_1w_3}}(w_{2}y_{1})B_{l,\chi,\xi^{w_2w_3}}
(w_{1}y_{2})B_{m,\chi,\xi^{w_1w_2}}(w_{3}y_{3}
)w_{2}^{l+m}w_{1}^{k+m}w_{3}^{k+l}\\
&=\sum_{k+l+m=n}^{} \binom{n}{k,l,m} B_{k,\chi
,\xi^{w_{1}w_{3}}}(w_{2}y_{1})B_{l,\chi,
\xi^{w_{1}w_{2}}}(w_{3}y_{2})B_{m,\chi,\xi^{w_{2}
w_{3}}}(w_{1}y_{3})w_{2}^{l+m}w_{3}^{k+m}
w_{1}^{k+l}\\
&=\sum_{k+l+m=n} \binom{n}{k,l,m} B_{k,\chi
,\xi^{w_{1}w_{2}}}(w_{3}y_{1})B_{l,\chi,
\xi^{w_{2}w_{3}}}(w_{1}y_{2})B_{m,\chi,\xi^{w
_{1}w_{3}}}(w_{2}y_{3})w_{3}^{l+m}w_{1}^{k+m}
w_{2}^{k+l}\\
&= \sum_{k+l+m=n} \binom{n}{k,l,m} B_{k,\chi
,\xi^{w_{1}w_{2}}}(w_{3}y_{1})B_{l,\chi,
\xi^{w_{1}w_{3}}}(w_{2}y_{2})B_{m,\chi,\xi^{w_{2}
w_{3}}}(w_{1}y_{3})w_{3}^{l+m}w_{2}^{k+m} w_{1}^{k+l}.
\end{split}
\end{equation}
}
\end{theorem}

\begin{theorem}\label{thm05}
Let $w_1,w_2,w_3$ be any positive integers, such that $r$ does not
divide $w_1w_2w_3$. Then we have: {\footnotesize
\begin{equation}\label{equ48}
\begin{split}
&\sum_{k+l+m=n} \binom{n}{k,l,m} B_{k,\chi,
\xi^{w_{2}w_{3}}(w_{1}y_{1})}B_{l,\chi,\xi^{w_{1}
w_{3}}}(w_{2}y_{2})S_{m}(dw_{3}-1;\chi
,\xi^{w_1w_2})w_{1}^{l+m}w_{2}^{k+m}w_{3}^{k+l}\\
&=\sum_{k+l+m=n}^{}
\binom{n}{k,l,m}B_{k,\chi,\xi^{w_{2}w_{3}}}(w_{1}y_{1})B_{l,\chi,\xi^{w_{1}w_{2}}}(w_{3}y_{2})S_{m}(dw_{2}-1;\chi,\xi^{w_{1}w_{3}})w_{1}^{l+m}w_{3}^{k+m}w_{2}^{k+l}
\\
&= \sum_{k+l+m=n}
\binom{n}{k,l,m}B_{k,\chi,\xi^{w_{1}w_{3}}}(w_{2}y_{1})B_{l,\chi,\xi^{w_{2}w_{3}}}(w_{1}y_{2})S_{m}(dw_{3}-1;\chi,\xi^{w_{1}w_{2}})w_{2}^{l+m}w_{1}^{k+m}w_{3}^{k+l}
\\
&= \sum_{k+l+m=n} \binom{n}{k,l,m}
B_{k,\chi,\xi^{w_{1}w_{3}}}(w_{2}y_{1})B_{l,\chi,\xi^{w_{1}w_{2}}}(w_{3}y_{2})S_{m}(dw_{1}-1;\chi,\xi^{w_{2}w_{3}})w_{2}^{l+m}w_{3}^{k+m}w_{1}^{k+l}
\\
&= \sum_{k+l+m=n} \binom{n}{k,l,m}
B_{k,\chi,\xi^{w_{1}w_{2}}}(w_{3}y_{1})B_{l,\chi,\xi^{w_{1}w_{3}}}(w_{2}y_{2})S_{m}(dw_{1}-1;\chi,\xi^{w_{2}w_{3}})w_{3}^{l+m}w_{2}^{k+m}w_{1}^{k+l}
\\
&= \sum_{k+l+m=n} \binom{n}{k,l,m}
B_{k,\chi,\xi^{w_{1}w_{2}}}(w_{3}y_{1})B_{l,\chi,\xi^{w_{2}w_{3}}}(w_{1}y_{2})S_{m}(dw_{2}-1;\chi,\xi^{w_{1}w_{3}})w_{3}^{l+m}w_{1}^{k+m}w_{2}^{k+l}.
\end{split}
\end{equation}
}
\end{theorem}

\begin{theorem}\label{thm06}
Let $w_1,w_2,w_3$ be any positive integers, such that $r$ does not
divide $w_1w_2w_3$. Then we have: {\footnotesize
\begin{equation}\label{equ49}
\begin{split}
& w_{1}^{n}\sum_{k=0}^{n} \binom{n}{k}B_{k,\chi,\xi
^{w_{1}w_{2}}}(w_{3}y_{1})\sum_{a=0}^{dw_{1}-1}
\chi(a)\xi^{aw_{2}w_{3}}B_{n-k,\chi,\xi^{w_{1}
w_{3}}}(w_{2}y_{2}+\frac{w_{2}}{w_{1}}a)w_{3}
^{n-k}w_{2}^{k}\\
&= w_{1}^{n}\sum_{k=0}^{n} \binom{n}{k}B_{k,\chi,
\xi^{w_{1}w_{3}}}(w_{2}y_{1})\sum_{a=0}^{dw_{1}
-1}\chi(a)\xi^{aw_{2}w_{3}}B_{n-k,\chi,\xi^{w
_{1}w_{2}}}(w_{3}y_{2}+\frac{w_{3}}{w_{1}}a)w
_{2}^{n-k}w_{3}^{k}\\
&= w_{2}^{n}\sum_{k=0}^{n} \binom{n}{k}B_{k,\chi,
\xi^{w_{1}w_{2}}}(w_{3}y_{1})\sum_{a=0}^{dw_{2}
-1}\chi(a)\xi^{aw_{1}w_{3}}B_{n-k,\chi,\xi^{w
_{2}w_{3}}}(w_{1}y_{2}+\frac{w_{1}}{w_{2}}a)w
_{3}^{n-k}w_{1}^{k} \\
&= w_{2}^{n}\sum_{k=0}^{n} \binom{n}{k}B_{k,\chi,
\xi^{w_{2}w_{3}}}(w_{1}y_{1})\sum_{a=0}^{dw_{2}
-1}\chi(a)\xi^{aw_{1}w_{3}}B_{n-k,\chi,\xi^{w
_{1}w_{2}}}(w_{3}y_{2}+\frac{w_{3}}{w_{2}}a)w
_{1}^{n-k}w_{3}^{k}\\
&= w_{3}^{n}\sum_{k=0}^{n} \binom{n}{k}B_{k,\chi,
\xi^{w_{1}w_{3}}}(w_{2}y_{1})\sum_{i=0}^{dw_{3}
-1}\chi(a)\xi^{aw_{1}w_{2}}B_{n-k,\chi,\xi^{w
_{2}w_{3}}}(w_{1}y_{2}+\frac{w_{1}}{w_{3}}a)w
_{2}^{n-k}w_{1}^{k}\\
&= w_{3}^{n}\sum_{k=0}^{n} \binom{n}{k}B_{k,\chi,
\xi^{w_{2}w_{3}}}(w_{1}y_{1})\sum_{a=0}^{dw_{3}
-1}\chi(a)\xi^{aw_{1}w_{2}}B_{n-k,\chi,\xi^{w
_{1}w_{3}}}(w_{2}y_{2}+\frac{w_{2}}{w_{3}}a)w _{1}^{n-k}w_{2}^{k}
\end{split}
\end{equation}
}
\end{theorem}

\begin{theorem}\label{thm07}
Let $w_1,w_2,w_3$ be any positive integers, such that $r$ does not
divide $w_1w_2w_3$. Then we have the following three symmetries in
$w_1,w_2,w_3$: {\tiny
\begin{align}
&\sum_{k+l+m=n} \binom{n}{k,l,m} B_{k,\chi,\xi^{w_{2}w_{3}}}(w_{1}y_{1})S_{l}(dw_{2}-1;\chi,\xi^{w_{1}w_{3}})S_{m}(dw_{3}-1;\chi,\xi^{w_{1}w_{2}})w_{1}^{l+m}w_{2}^{k+m}w_{3}^{k+l}\label{equ50}\\
&= \sum_{k+l+m=n} \binom{n}{k,l,m}B_{k,\chi,\xi^{w_{1}w_{3}}}(w_{2}y_{1})S_{l}(dw_{3}-1;
\chi,\xi^{w_{1}w_{2}})S_{m}(dw_{1}-1;\chi,\xi^{w_{2}w_{3}})w_{2}^{l+m}w_{3}^{k+m}w_{1}^{k+l}\label{equ51}\\
&= \sum_{k+l+m=n} \binom{n}{k,l,m}B_{k,\chi,
\xi^{w_{1}w_{2}}}(w_{3}y_{1})S_{l}(dw_{1}-1;
\chi,\xi^{w_{2}w_{3}})S_{m}(dw_{2}-1;\chi,\xi
^{w_{1}w_{3}})w_{3}^{l+m}w_{1}^{k+m}w_{2}^{k+l}.\label{equ52}
\end{align}
}
\end{theorem}

\begin{theorem}\label{thm08}
Let $w_1,w_2,w_3$ be any positive integers, such that $r$ does not
divide $w_1w_2w_3$. Then we have: {\footnotesize
\begin{equation}\label{equ53}
\begin{split}
& w_{1}^{n}\sum_{k=0}^{n} \binom{n}{k}\sum_{a=0}^{dw_{1}
-1}\chi(a)\xi^{aw_{2}w_{3}}B_{k,\chi,\xi^{aw
_{1}w_{3}}}(w_{2}y_{1}+\frac{w_{2}}{w_{1}}a)S
_{n-k}(dw_{3}-1;\chi,\xi^{w_{1}w_{2}})w_{2}
^{n-k}w_{3}^{k}\\
&= w_{1}^{n}\sum_{k=0}^{n} \binom{n}{k}\sum_{a=0}^{dw
_{1}-1}\chi(a)\xi^{aw_{2}w_{3}}B_{k,\chi,\xi^{w
_{1}w_{2}}}(w_{3}y_{1}+\frac{w_{3}}{w_{1}}a)S
_{n-k}(dw_{2}-1;\chi,\xi^{w_{1}w_{3}})w_{3}
^{n-k}w_{2}^{k}\\
&= w_{2}^{n}\sum_{k=0}^{n} \binom{n}{k}\sum_{a=0}^{dw
_{2}-1}\chi(a)\xi^{aw_{1}w_{3}}B_{k,\chi,\xi^{w
_{2}w_{3}}}(w_{1}y_{1}+\frac{w_{1}}{w_{2}}a)S
_{n-k}(dw_{3}-1;\chi,\xi^{w_{1}w_{2}})w_{1}
^{n-k}w_{3}^{k}\\
&= w_{2}^{n}\sum_{k=0}^{n} \binom{n}{k} \sum_{a=0}^{dw
_{2}-1}\chi(a)\xi^{aw_{1}w_{3}}B_{k,\chi,\xi^{w
_{1}w_{2}}}(w_{3}y_{1}+\frac{w_{3}}{w_{2}}a)S
_{n-k}(dw_{1}-1;\chi,x^{w_{2}w_{3}})w_{3}^{n-k}
w_{1}^{k}\\
&= w_{3}^{n}\sum_{k=0}^{n} \binom{n}{k}\sum_{a=0}^{dw
_{3}-1}\chi(a)\xi^{aw_1w_2}B_{k,\chi,\xi^{w_2w_3}
}(w_{1}y_{1}+\frac{w_{1}}{w_{3}}a)S_{n-k}(dw_{2}
-1;\chi,\xi^{w_1w_3})w_{1}^{n-k}w_{2}^{k}\\
&= w_{3}^{n}\sum_{k=0}^{n} \binom{n}{k}\sum_{a=0}^{dw
_{3}-1}\chi(a)\xi^{aw_{1}w_{2}}B_{k,\chi,\xi
^{w_{1}w_{3}}}(w_{2}y_{1}+\frac{w_{2}}{w_{3}}
a)S_{n-k}(dw_{1}-1;\chi,\xi^{w_{2}w_{3}})w_{2}^{n-k}w_{1}^{k}.
\end{split}
\end{equation}
}
\end{theorem}

\begin{theorem}\label{thm09}
Let $w_1,w_2,w_3$ be any positive integers, such that $r$ does not
divide $w_1w_2w_3$. Then we have the following three symmetries in
$w_1,w_2,w_3$: {\small
\begin{equation}\label{equ54}
\begin{split}
&(w_1w_2)^{n}\sum_{a=0}^{dw_1-1}\sum_{b=0}^{dw_2-1}\chi(ab)\xi^{w_3(aw_2+bw_1)}B_{n,\chi,\xi^{w_1w_2}}
(w_3y_1+\frac{w_3}{w_1}a+\frac{w_3}{w_2}b
)\\
&= (w_{2}w_{3})^{n}\sum_{a=0}^{dw_{2}-1}\sum_{b=0}^{dw_{3}-1}
\chi(ab)\xi^{w_{1}(aw_{3}+bw_{2})}B_{n,\chi,\xi^{w_{2}w_{3}}}
(w_{1}y_{1}+\frac{w_{1}}{w_{2}}a+\frac{w_{1}}{w_{3}}
b) \\
&= (w_3w_1)^{n}\sum_{a=0}^{dw_3-1}\sum_{b=0}^{dw_1-1}
\chi(ab)\xi^{w_2(aw_1+bw_3)}B_{n,\chi,\xi^{w_1w_3}}(w_2y_1+
\frac{w_2}{w_3}a+\frac{w_2}{w_1}b).
\end{split}
\end{equation}
}
\end{theorem}

\begin{theorem}\label{thm10}
Let $w_1,w_2,w_3$ be any positive integers, such that $r$ does not
divide $w_1,w_2,w_3$. Then we have the following two symmetries in
symmetries in $w_1,w_2,w_3$:

{\small
\begin{equation}\label{equ55}
\begin{split}
& \sum_{k+l+m=n} \binom{n}{k,l,m}  B_{k,\chi,\xi^{w _{3}}}(w_{1}y)
 B_{l,\chi,\xi^{w_{1}}}(w_{2}y) B_{m,\chi,\xi^{w_{2}}}(w_{3}y)w_3^kw_1^lw_2^m \\
&\quad = \sum_{k+l+m=n} \binom{n}{k,l,m}  B_{k,\chi,\xi ^{w_{2}}}(w_{1}y)
 B_{l,\chi,\xi^{w_{1}}}(w_{3} y) B_{m,\chi,\xi^{w_{3}}}(w_{2}y)w_2^k w_1^lw_3^m.
\end{split}
\end{equation}
}

\end{theorem}

\begin{theorem}\label{thm11}
Let $w_1,w_2,w_3$ be any positive integers, such that $r$ does not
divide $w_2w_3,w_1w_3,w_1w_2$. Then we have the following two
symmetries in $w_1,w_2,w_3$: {\footnotesize
\begin{align}
& \sum_{k+l+m=n}  \binom{n}{k,l,m}
S_{k}(dw_{1}-1;\chi,\xi^{w_{3}}) S_{l}(dw_{2}-1;\chi,\xi^{w_{1}}) S_{m}(dw_{3}-1;\chi,\xi^{w_{2}})w_3^{k} w_1^{l}w_2^{m}\label{equ56}\\
& =\sum_{k+l+m=n}  \binom{n}{k,l,m}
S_{k}(dw_{1}-1;\chi,\xi^{w_2})S_{l} (dw_{3}-1;\chi,\xi^{w_1})
S_{m}(dw_{2}-1;\chi,\xi^{w_3})
w_2^{k}w_1^{l}w_3^{m}.\label{equ57}
\end{align}
}
\end{theorem}

%%%%%%%%%%%%%%%%%%%%%%%%%%%%%%%%%%%%%%%%%%%%%%%%%%%%%%%%%%%%%%%%%%%%%%%%%%%%%%%%%%%%%%%%%%
%%%%%%%%%%%%%%%%%%%%%%%%%%%%%%%%%%%%%%%%%%%%%%%%%%%%%%%%%%%%%%%%%%%%%%%%%%%%%%%%%%%%%%%%%%
%%%%%%%%%%%%%%%%%%%%%%%%%%%%%%%%%%%%%%%%%%%%%%%%%%%%%%%%%%%%%%%%%%%%%%%%%%%%%%%%%%%%%%%%%%

%%%%%%%%%%%%%%%%%%%%%%%%%%%%%%%%%%%%%%%%%%%%%%%%%%%%%%%%%%%%%%%%%%%%%%%%%%%%%%%%%%%%%%%%%%
%%%%%%%%%%%%%%%%%%%%%%%%%%%%%%%%%%%%%%%%%%%%%%%%%%%%%%%%%%%%%%%%%%%%%%%%%%%%%%%%%%%%%%%%%%
%%%%%%%%%%%%%%%%%%%%%%%%%%%%%%%%%%%%%%%%%%%%%%%%%%%%%%%%%%%%%%%%%%%%%%%%%%%%%%%%%%%%%%%%%%

\end{document}